\documentstyle[12pt]{article}

\input{mssymb}

\newtheorem{lemma}{Lemma}
      
\newtheorem{thm}{Theorem}

\newcommand{\re}{\restriction}
\newcommand{\be}{\begin{enumerate}}
\newcommand{\ee}{\end{enumerate}}
\newcommand{\bq}{\begin{quote}}
\newcommand{\eq}{\end{quote}}

\newcommand{\seq}{\subseteq}
\newcommand{\calC}{{\cal C}_{\omega}}

\newpage

\begin{document}
\author{A. Louveau, S. Shelah, and B. Velickovic}
\title{Borel partitions of infinite sequences of reals}
\date{}
\maketitle

\section*{\bf Introduction}

The starting point of our work is a Ramsey-type theorem of
Galvin (unpublished) which asserts that if the unordered pairs
of reals are partitioned into finitely many Borel classes (or even
classes which have the property of Baire) then there is a perfect
set $P$ such that all pairs from $P$ lie in the same class.
The obvious generalization to $n$-tuples for $n\geq 3$
is false. For example, look at the coloring of triples where
a triple $\{x,y,z\}$ with $x<y<z$ is colored {\em red} provided
that $y-x < z-y$ and {\em blue} otherwise. Then any perfect 
set will contain triples of both colors. Galvin conjectured
that this is the only bad thing that can happen. It will 
be simpler to state this if we identify the reals with
$2^{\omega}$ ordered by the lexicographical ordering 
and define for distinct $x,y\in 2^{\omega}$  $\Delta (x,y)$
to be the least $n$ such that $x(n)\neq y(n)$.
Let the {\em type} 
of an increasing $n$-tuple of reals $\{x_0,\ldots x_{n-1}\}_<$ be
the ordering $\prec$ on $\{0,\ldots ,n-2\}$ defined by 
$i\prec j\ \mbox{iff}\ \Delta (x_i,x_{i+1})<\Delta (x_j,x_{j+1}).$
Galvin proved that for any Borel coloring 
of triples of reals there is a perfect set $P$ such
that the color of any triple from $P$ depends only on its type
and conjectured that an analogous result is true for any $n$.
This conjecture  has been proved by Blass ([Bl]). As a corollary
it follows that if the unordered $n$-tuples of reals are
colored into finitely many Borel classes there is a perfect
set $P$ such that the $n$-tuples from $P$ meet at most $(n-1)!$ classes.
The key ingredient in the proof is the well-known
Halpern-La\"{u}chli theorem ([HL]) on partitions of products
of finitely many tree. 
In this paper we consider extensions of this result to partitions
of infinite increasing sequences of reals.
Define a type of an increasing sequence of reals as before and
say that such a sequence $\{x_n:n<\omega \}$ is {\em strongly increasing}
if its type is the standard ordering on $\omega$, i.e.
if  $\Delta (x_n,x_{n+1})<\Delta (x_{m},x_{m+1})$ whenever $n<m$.
We show, for example, that for any Borel or even analytic
partition of all  increasing sequences
of reals there is a perfect set $P$ such that all strongly
increasing sequences from $P$ lie in the same class.
In fact, for any finite set ${\cal C}$ of types
there is a perfect set $P$ such that
for any type in ${\cal C}$ all increasing sequence
from $P$ of that type have the same color. It should be pointed out that 
the same statement is false if ${\cal C}$ is an infinite set of types.

Our result stands in the same relation to Blass' theorem
as the Galvin-Prikry theorem ([GP]) to the ordinary
Ramsey's theorem and the proof again relies heavily on
the Halpern-La\"{u}chi theorem. There are known several extensions
of the Halpern-La\"{u}chli theorem that are relevant to this work.
Milliken ([Mi]) considered partitions of nicely embedded infinite subtrees
of a perfect tree and obtained a partition result in the spirit
of Galvin-Prikry however in a different direction from ours,
and Laver ([La]) proved a version of this theorem for products
of infinitely many perfect trees.  

The paper is organizes as follows. In \S 1 we introduce some notation
and present some results on perfect trees which we will need later.
In  \S 2 we reduce the main theorem to two lemmas which are 
then proved in \S \S 3 and 4. We shall present our result using
the terminology of forcing. If ${\cal P}$ is a forcing notion
we let, as usual, $RO({\cal P})$ denote the regular open algebra of ${\cal P}$,
i.e. a complete Boolean algebra in which ${\cal P}$ is densely embedded.
If ${\bf b}$ is a Boolean value in $RO({\cal P})$
and $p\in {\cal P}$ we shall say that $p$ {\em decides} ${\bf b}$
if either $p\leq {\bf b}$ or $p\leq {\bf 1-b}$.
For all undefined terminology of forcing see, for example, [Ku].

\section{\bf Basic properties of perfect trees}

\noindent {\bf Perfect trees} Let  $2^{< \omega}$\ denote the set of 
all finite $\{0,1\}$-sequences ordered by extension. $T \subseteq 2^{< \omega}$
is called a {\em perfect tree} if it is an initial segment of $2^{< \omega}$
and every element of $T$ has two incomparable extensions in $T$.
Let ${\cal P}$\ denote the poset  of all perfect trees partially ordered
by inclusion. Thus ${\cal P}$\ is the well-known Sacks forcing ([Sa]).
For a subset $C$ of $T$ let
$T_C$ be the set of all nodes in $T$ which are comparable to an element of $C$.
If $\{s\}$ is a singleton  we shall simply write $T_s$ instead of $T_{\{s\}}$.
For a tree $T$ let $T(n)$ denote the $n$-th level of $T$, i.e. the set
of all $s\in T$ which have exactly $n$ predecessors. We say that a node
$s$ in $T$ is {\em splitting} if it has two immediate extensions. 
Given integers $m\leq k$ let us say that a set $D$ is $(m,k)$-{\em dense} in $T$
provided $D$ is contained in $T(k)$ and every node in $T(m)$ has an extension in $D$.
Given trees $T_0,\ldots T_{d-1}$ and a subset $A$ of $\omega$ let
$$
\otimes ^A _{i<d}T_i = \bigcup _{n\in A}\otimes _{i<d}T_i(n).
$$
If $A$ is $\omega$ we usually omit it.
 We are now ready to state a version of the Halpern-La\"{u}chli theorem ([HL]).
\begin{thm}[{[}HL{]}]{ For every integer $d<\omega$ given  perfect trees $T_i$,
for $i<d$,  and a partition 
$$
\otimes_{i<d}T_i = K_0 \cup K_1
$$
for every  infinite subset $A$ of $\omega$ 
there are $(x_0,\ldots ,x_{d-1})\in \otimes _{i<d}T_i$ and $\epsilon \in \{0,1\}$ such 
that for every $m$ there is $k\in A$ and sets $D_i$, for $i<d$, such that
$D_i$ is $(m,k)$-dense in $T_i$ and  $\otimes _{i<d}D_i\seq K_{\epsilon}$.} 
\end{thm}

\medskip

\noindent {\bf The amoeba forcing ${\cal A}({\cal P})$} To the poset 
 ${\cal P}$\ we associate the {\em amoeba poset}
${\cal A}({\cal P})$. Elements of ${\cal A}({\cal P})$\ are pairs
$(T,n)$, where $T \in {\cal P}$ and  $n \in \omega$.
Say that $(T,n) \leq (S,m)$\ iff $T \leq S, n \geq m$, and 
$T \re (m+1) = S \re (m+1).$ If in addition $n=m$ we shall say that
$(T,n)$ is a {\em pure extension} of $(S,m)$.
If $G$ an ${\cal A}({\cal P})$-generic filter over a model of set theory let
$$
T(G) = \bigcup \{ T \re (n+1) : (T,n) \in G \}.
$$
Then, by genericity,  $T(G)$\ is a perfect tree and is  called the 
${\cal A}({\cal P})$-{\em generic tree} derived from $G$.

\medskip 

\noindent {\bf Combs} An $n$-{\em comb} $C$ is a tree
such that there is some strongly increasing sequence of reals $\{ x_i:i<n\}$ and some
$m>\Delta (x_{n-2},x_{n-1})$ such that $C$ is the set of all initial
 segments of length $<m$ of members of this sequence. 
An {\em infinite comb} is a tree such that there is  some strongly increasing
sequence $\{ x_n:n<\omega \}$ such that  $C$ is the set of all finite 
initial segments of members of this sequence.
 Clearly there is a 1-1 correspondence
between infinite combs and strongly increasing sequences and we shall 
in fact state our theorem in terms of infinite combs. For a tree
$T$ if $n<\omega$ is such that $T\re (n+1)$ is  a comb let
 ${\cal C}_{\omega}(T,n)$ denote the set of all infinite combs
contained in $T$ and extending $T\re (n+1)$. 
Let $\calC (T)= \calC (T,0)$.  Note that ${\cal C}_{\omega}(T)$ has 
a natural topology as a subspace of
${\cal P}(T)$ with the Tychonoff topology. 
Thus we can speak about Borel, analytic, etc.
subsets of ${\cal C}_{\omega}(T)$. 

\medskip

\noindent {\bf The comb forcing ${\cal C}$} Let ${\cal C}$ be the subposet of
${\cal A}({\cal P})$ consisting of all pairs $(T,n)$ such that
$T\re (n+1)$ is a  comb, with the induced ordering. Let us say that 
$(T,n)$ has width $d$ if $T\re (n+1)$ is a $d$-comb.
The notion of pure extension is defined as in the case of 
${\cal A}({\cal P})$.  If $(R,m) \leq (T,n)$ and if these two conditions 
have the same width then we say that $(R,m)$ is a {\em width preserving}
extension of $(T,n)$. Note that in this case $(R,n)$ is a pure 
extension of $(T,n)$ which is equivalent in terms of forcing with $(R,m)$.
Clearly, if $G$ is a ${\cal C}$-generic filter over some model of
set theory  the set
$$
C(G)=\bigcup \{ T\re (n+1): (T,n) \in G \}
$$
is a infinite comb, we call it the {\em generic comb} derived from $G$.

\medskip

\section{\bf The main theorem}
The main result of this paper is the following partition theorem.

\begin{thm}{For every partition
$$
\calC (2^{<\omega})=K_0 \cup K_1
$$
where $K_0$ is analytic and $K_1$ co-analytic there is a perfect
tree $T$ and $i\in \{ 0,1\}$ such that $\calC (T)\seq K_i$.}
\end{thm}

The proof of the theorem will consist of two lemmas which combined
yield the desired result.

\begin{lemma}{Let ${\bf b}$ be a Boolean value in $RO({\cal C})$
and let $(S,n)\in {\cal C}$. Then there is a pure extension
$(T,n)$ of $(S,n)$ which decides ${\bf b}$.}
\end{lemma}

\begin{lemma}{Let $T$ be an ${\cal A}({\cal P})$-generic tree 
over a model of set theory $M$. Then every infinite comb contained
in $T$ is ${\cal C}$-generic over $M$.}
\end{lemma}

Given these two lemmas it is quite easy to prove the theorem.
Take a countable transitive model $M$ of ZFC$^-$ containing 
the codes of $K_0$ and $K_1$. Consider forcing with ${\cal C}$ 
as defined in $M$. Note that if $C$ is a generic comb
the statement whether $C$ belongs to $K_0$ is absolute between
$M[C]$ and $V$. Let ${\bf b}$ be the Boolean value that this 
statement is true in $M[C]$. Then it follows from Lemma 1 that there is a pure
extension $(S,0)$ of the maximal condition which decides 
${\bf b}$, let us say, for concreteness,  that it forces ${\bf b}$.
Now consider forcing over $M$ with ${\cal A}({\cal P})$ and 
take a generic filter $G$ over $M$ which contains $(S,0)$.
Let $T$ be the generic tree derived from $G$. Then by Lemma 2
every infinite comb contained in $T$ is ${\cal C}$-generic 
over $M$ and, since it is contained in $S$ as well,  it 
follows that it is in $K_0$. Thus $T$ is the homogeneous 
tree we seek.
In the next two sections we prove Lemmas 1 and 2 and thus complete 
the proof.

\section{\bf Proof of Lemma 1}

Unless otherwise stated in this section we work with 
 the forcing notion ${\cal C}$ introduced in \S 1.
Given a Boolean value ${\bf b}$ in the completion algebra $RO({\cal C})$
let us say that a condition $(T,n)$ {\em accepts} ${\bf b}$
if $(T,n) \leq{\bf b}$ and that it {\em rejects} ${\bf b}$ if
$(T,n)\leq {\bf 1-b}$. We shall need the following auxiliary lemma.

\begin{lemma}{Let $(S,n)$ be a condition in ${\cal C}$ of width $d$  and let
${\bf b} \in RO({\cal C})$ be a Boolean value.
Then there is a pure extension $(T,n)$ of $(S,n)$ such
that either $(T,n)$ accepts ${\bf b}$ or
 no extension of $(T,n)$ of width $d+1$ accepts ${\bf b}$.}
\end{lemma}

\noindent PROOF: Let $\{t_0,\ldots , t_{d-1}\}_<$ be the increasing enumeration
of $S(n)$ in the lexicographical ordering. We first find an infinite
set $A$ and a  perfect subtree $S^*$ of $S$ such that for any $m\in A$
and $z_0,\ldots ,z_d\in S^*(m)$ such that $z_i\geq t_i$ for $i<d$
and $z_d\geq t_{d-1}$, letting $Z=\{z_i:i\leq d\}$, if there is
a pure extension of $(S^*_Z,m)$ deciding ${\bf b}$ then already 
$(S^*_Z,m)$ decides ${\bf b}$. This can be done by a standard fusion argument.
Moreover we can arrange that between any two consecutive levels
in $A$ there is at most one splitting node. We now define a coloring:
$$
\otimes _{i<d} S^*_{t_i}=K_0\cup K_1\cup K_2
$$
as follows. Given $(x_0,\ldots ,x_{d-1})\in \otimes _{i<d}S^*_{t_i}$
let $m\in A$ be the least such that $x_{d-1}$ has two extensions
$z_{d-1}$ and $z_d$ in $S^*(m)$. For $i<d-1$ let
$z_i$ be the lexicographically least extension of $x_i$ in $S^*(m)$.
Let $Z=\{z_i:i\leq d\}$ and put $(x_0,\ldots ,x_{d-1})$
in $K_0$ if $(S^*_Z,m)$ accepts ${\bf b}$, 
in $K_1$ if it rejects ${\bf b}$, and in $K_2$ otherwise.
By the Halpern-La\"{u}chli theorem we can find
$(x_0,\ldots ,x_{d-1})\in \otimes _{i<d}S^*_{t_i}$ and 
$\epsilon \in \{0,1,2\}$ such that for every $m$ there is $k\in A$
and sets $D_i$, for $i<d$, such that $D_i$ is $(m,k)$-dense in
$S^*_{x_i}$ and $\otimes _{i<d} D_i\seq K_{\epsilon}$.
We may assume that $(x_0,\ldots ,x_{d-1})\in K_{\epsilon}$, as well.

We now build an increasing sequence $(b_k)_{k<\omega}$
of elements of $A$ and a perfect subtree $T$ of $S^*$  which
will have one splitting node on levels  between 
$b_k$ and $b_{k+1}$.
To begin let $b_0$ be the level of the $x_i$ and let 
$T(b_0)=\{x_i:i<d\}$. This uniquely determines $T\re (b_0+1)$ as the set of all
initial segments of elements of $T(b_0)$.
Suppose now we have defined $b_k$ and $T\re (b_k+1)$.
We choose one node $y$ in $T(b_k)$ and we will arrange so that
the only splitting node of $T$ on levels between $b_k$ and $b_{k+1}$ is above $y$.
Let $m$ be the least level which is in $A$ and such that $y$ has two
extensions, say $y^{\prime}$ and $y^{\prime \prime}$  in $S^*(m)$.
Now find some $b \in A$ and sets $D_i$, for $i<d$ such that 
$D_i$ is $(m,b)$-dense in $S^*_{x_i}$ and such that 
$\otimes _{i<d}D_i\seq K_{\epsilon}$. Set $b_{k+1}=b$ and let
$D=\bigcup _{i<d}D_i$. For each element in 
$T(b_k)\cup \{y^{\prime},y^{\prime \prime}\}$ pick a lexicographically
least point in $D$ above it. Let $T(b_{k+1})$ be the set of 
points thus chosen. This uniquely defines $T\re (b_{k+1}+1)$.
During our construction we arrange the choice of the points $y$ in such
a way that the final tree $T$ is perfect. Let $B=\{b_k: k<\omega \}$.
It follows that $\otimes _{i<d}^BT_{t_i} \seq K_{\epsilon}$.

We now show that $(T,n)$ is the required condition.
First note that if $(R,l)$ is any extension of $(T,n)$ then there 
there is $m\in A$ such that  $R$ has no splitting nodes	
on levels between $l$ and $m$ and  hence $(R,l)$ and $(R,m)$ are equivalent condition.
Suppose now  that some condition of width $d+1$ below $(T,n)$ 
accepts ${\bf b}$ and let $(R,m)$ be such a condition with $m$ minimal such that $m\in A$.
Let $Z=R(m)=\{z_0,\ldots ,z_d\}_<$ be the increasing enumeration in
the lexicographical order and let $k$ be the largest such that
$b_k <m$. Then since on levels between $b_k$ and $b_{k+1}$ 
there is at most one  splitting node it follows that
$R(b_k)$ has size $d$. Let $R(b_k)=\{y_0,\ldots ,y_{d-1}\}_<$ be
the increasing enumeration. By the construction of $T$ 
it follows that $y_{d-1}$ was the point chosen at stage $k$,
that $z_{d-1}$ and $z_d$ are the only extensions of $y_{d-1}$
in $T$ on level $m$, and that $z_i$ is the lexicographically least extension
of $y_i$ in  $S^*(m)$ for $i<d-1$.
Thus $(y_0,\ldots ,y_{d-1})$ is colored according to whether 
$(S_Z^*,m)$ accepts ${\bf b}$, rejects ${\bf b}$, or cannot decide.
Since $(R,m)$ is a pure extension of $(S^*_Z,m)$ which
accepts ${\bf b}$ and $m\in A$  by the property of $S^*$ 
it follows that $(S^*_Z,m)$ also accepts ${\bf b}$ 
and thus $(y_0,\ldots ,y_{d-1})\in K_0$. Hence we must have 
$\epsilon =0$. 

Now since then $\otimes _{i<d}^B T_{t_i} \seq K_0$,
a similar analysis shows that any other extension of 
$(T,n)$ of width $d+1$ accepts ${\bf b}$. But then it follows that
$(T,n)$ also accepts ${\bf b}$. 
\hspace{.25in} $\Box$
\medskip

\noindent PROOF OF LEMMA 1: Let $(S,n)$ be a condition  in ${\cal C}$
and let ${\bf b}$ be a Boolean value. Assume that there is no
pure extension of $(S,n)$ which accepts ${\bf b}$. We find 
a pure extension $(T,n)$ of $(S,n)$ which rejects ${\bf b}$.
We shall build the tree $T$ by a fusion argument. 
Along the way we shall construct a decreasing sequence 
$(T^{(0)},a_0)\geq (T^{(1)},a_1) \geq \ldots$
of conditions in ${\cal A}({\cal P})$.

To begin let  $(T^{(0)},a_0)=(S,n)$.
Suppose now $(T^{(k)},a_k)$, has been defined.
Let $\{Z_i:i<l\}$ be an enumeration of all subsets $Z$ of 
$T^{(k)}(a_k)$ which generate a comb extending $S\re (n+1)$.
The inductive assumption is that for each such $Z$ the condition 
$(T^{(k)}_Z,a_k)$ does not have a pure extension accepting ${\bf b}$.
To avoid excessive notation let $R$ be a variable denoting a perfect
subtree of $T^{(k)}$. We initially set $R$ to be equal $T^{(k)}$ and
then trim it down in $l$ steps as follows. 
At step $i$ consider $Z_i$. 
Since $(R_{Z_i},a_k)$ is a pure extension of $(T^{(k)}_{Z_i},a_k)$
from the inductive assumption it follows that it does not have a pure
extension accepting ${\bf b}$. 
If the size of $Z_i$ is $d_i$ then by Lemma 3 there is a pure extension
$(Q,a_k)$ of $(R_{Z_i},a_k)$ such that
no extension of $(Q,a_k)$ of width $d_i+1$ accepts ${\bf b}$.
We now shrink $R$ as follows. For every $s\in Z_i$ replace $R_s$ by $Q_s$
and for $s\in T^{(k)}(a_k)\setminus Z_i$ keep $R_s$ the same.
After all the $l$ steps have been completed 
pick a node $y$ in $T^{(k)}(a_k)$.
Let $a_{k+1}$ be the least $a$  such that $y$ has two extensions in $R(a)$.
Keep those two extensions of $y$ and for every other node in $T^{(k)}(a_k)$ pick
exactly one extension on level $a_{k+1}$. Let then
 $T^{(k+1)}$ be the set of all nodes of $R$ comparable to one of  these nodes.
If now $Z$ is any subset of $T^{(k+1)}(a_{k+1})$ which generates a comb
extending $S\re (n+1)$ we claim that there is no pure extension of 
$(T^{(k+1)},a_{k+1})$ accepting ${\bf b}$.
Notice that the set of all predecessors of members of $Z$ on level $a_k$
is listed as one of the $Z_i$. Since between levels $a_k$ and $a_{k+1}$ 
there is at most one splitting of $T^{(k+1)}$ it follows that 
card$(Z)\leq d_i+1$.  If the size of $Z$ is $d_i$ then 
every pure extension of $(T^{(k+1)},a_{k+1})$ is equivalent to 
a pure extension of  $(T^{(k)},a_k)$, but by the inductive hypothesis 
such a condition cannot accept ${\bf b}$. On the other hand
if the size of $Z$ is $d_i+1$ at stage $i$ of the construction of $T^{(k+1)}$
we have ensured that no such condition accepts ${\bf b}$.
This shows that the inductive hypothesis is preserved.

Finally let $T =\bigcap T^{(k)}$. 
Throughout the construction we make the choice of the points $y$ above which
we keep a splitting node carefully to ensure that the final tree $T$ is perfect. 
It follows that no condition $(R,m)$ extending $(T,n)$  accepts ${\bf b}$
and hence $(T,n)$ rejects ${\bf b}$, as desired. \hspace{.25in} $\Box$

\section{\bf Proof of Lemma 2}

In the proof of Lemma 2 we need the following lemma whose proof is
almost identical to the proof of Lemma 3 and is thus omitted.

\begin{lemma}{ Let $(S,n)\in {\cal C}$ be a condition of width $d$ and
let $U$ be a set of infinite combs. Then there is pure extension
$(T,n)$ of $(S,n)$ such that either $\calC (T,n)$ is contained in $U$ or 
there is no extension $(R,m)$ of $(T,n)$ of width $d+1$ 
such that $\calC (R,m)$ is contained in $U$.}
\end{lemma}

Now note that to  complete the proof of Lemma 2 and Theorem 2 it suffices
to prove the following. 

\begin{lemma}{Let $(S,n)$  be a condition in ${\cal A}({\cal P})$ and let $D$ be
a dense open subset of ${\cal C}$. Then there is a pure extension
$(T,n)$ of $(S,n)$ such that for every infinite comb $C$ in $\calC (T)$
there is $m$ such that $(T_{C(m)},m)\in D$.} 
\end{lemma}

\noindent PROOF: We first show  that if $(S,n)\in {\cal C}$ there is a pure
extension $(T,n)$ of $(S,n)$ such that for every $C\in \calC (T,n)$ there
is $m\geq n$ such that $(T_{C(m)},m)\in D$. 
To begin find an infinite subset $A$ of $\omega$ and a pure extension
$(S^*,n)$ of $(S,n)$ such that for every $m\in A$ and every subset 
$Z$ of $S^*(m)$ which generates a comb extending $S\re (n+1)$ 
if there is a pure extension of $(S^*_Z,m)$ which is in $D$ 
then already $(S^*_Z,m)$ is in $D$. Let then
$$
U=\{ C\in \calC (S^*,n) :\ \mbox{there is $m$ such that}\ (S^*_{C(m)},m) \in D\}
$$
Assume now towards contradiction that there is no pure extension 
$(T,n)$ of $(S^*,n)$ such that $\calC (T,n)$ is contained in $	U$. As in the proof
of Lemma 1 we build a decreasing sequence 
$(T^{(0)},a_0)\geq (T^{(1)},a_1)\geq \ldots$ of conditions in ${\cal A}({\cal P})$.
 To begin set $(T^{(0)},a_0)=(S^*,n)$.
Suppose now $(T^{(k)},a_k)$ has been defined. Our inductive assumption 
is that for any subset $Z$ of $T^{(k)}(a_k)$ which generates a comb
extending $S\re (n+1)$ there is no pure extension $(Q,a_k)$ of 
$(T^{(k)}_Z,a_k)$ such that $\calC (Q,a_k)$  is contained in $U$.
Let $\{Z_i:i<l\}$ be an enumeration of all such $Z$.
To avoid excessive notation let, as before, $R$ be a variable 
denoting a perfect subtree of $T^{(k)}$.
To begin set $R$ equal to $T^{(k)}$. We then successively trim down $R$ 
in $l$ steps as follows. Suppose that step $i$ has been completed.
Since $(R_{Z_i},a_k)$ is a pure extension of $(T^{(k)}_{Z_i},a_k)$, by the inductive
hypothesis it has  no pure extension $(Q,a_k)$ such that $\calC (Q,a_k)$
is contained in $U$.
Let the size $Z_i$ be $d_i$. 
Then by Lemma 4 there is a pure extension $(Q,a_k)$ of $(R_{Z_i},a_k)$
in ${\cal C}$  such that if $(Q^*,m)$ is an  extension $(Q,a_k)$ in  ${\cal C}$ 
of width $d_i+1$ then  $\calC (Q^*,m)$ is not contained in $U$.
Now trim down $R$ as follows. For nodes $s$ in $Z_i$ replace $R_s$ 
by $Q_s$ and for nodes $s$ in $T^{(k)}(a_k)\setminus Z_i$ keep $R_s$ the same.
Finally when all the stages are completed and we have taken care of
all the $Z_i$ we choose a node $y$ in $T^{(k)}(a_k)$ and let 
$a_{k+1}$ be the least member of $A$ above $a_k$ such that 
$y$ has two successors in $R(a_{k+1})$. Then $T^{(k+1)}$ is
obtained from $R$ by keeping those two successors of $y$ and by keeping 
for every other node in $T^{(k)}(a_k)$ one successors and throwing away 
the remaining ones. Then $T^{(k+1)}$ is set to be the set of
all nodes of the final $R$ comparable to one of the chosen points.
Note that in this way we arrange that for every subset $Z$ of 
$T^{(k+1)}(a_{k+1})$ which generates a comb extending $S\re (n+1)$ 
the set of all  predecessors of members of $Z$ on level $a_k$ is listed as
one of the $Z_i$ and since between $a_k$ and $a_{k+1}$ there is 
at most one splitting node it follows that card$(Z)\leq d_i+1$.
Thus it follows that if $(Q,a_{k+1})$ is a pure extension 
of  $(T^{(k+1)}_Z,a_{k+1})$ then   $\calC (Q,a_{k+1})\setminus U\neq \emptyset$.

In then end we let $T= \bigcap T_k$. We make the choice of the nodes
$y$ above we choose a splitting at each stage judiciously so that the final
tree  $T$ is perfect. It follows that if $(R,m)$ is any extension 
of $(T,n)$ in ${\cal C}$ then $\calC (R,m)\setminus U \neq \emptyset$. 
Now since $D$ is dense open we can find $k$ and a condition
$(R,a_k)\in D$ extending $(T,n)$. Let $Z=R(a_k)$. By the property of 
$S^*$ it follows that $(S^*_Z,a_k)$ is also in $D$. But then 
$\calC (S^*_Z,a_k)\seq U$, a contradiction. 

Now to deal with the general case assume that only $(S,n)\in {\cal A}({\cal P})$.
We then proceed as in the successor stage of the previous case.
We enumerate all subset $Z$ of $S(n)$ which generate a comb
as $\{ Z_i:i<l\}$.
Let, as before, $R$ be a variable denoting a perfect subtree of 
$S$. To begin set $R$ to be equal to $S$. 
We then trim down $R$ successively in $l$ stages. At stage
$i$ look at $Z_i$ and apply the special case of the lemma to 
find a pure extension $(Q,n)$ of $(R_{Z_i},n)$ such that
for every infinite comb $C$ extending $Q\re (n+1)$ there is
$m\geq n$ such that $(Q_{C(m)},m)\in D$.
Trim down $R$ by replacing $R_s$ by $Q_s$ 
for every node $s\in Z_i$ and keeping $R_s$ the same for 
evert $s\in S(n)\setminus Z_i$. 
We let $T$ be equal to $R$ after all the stages have been completed.
It follows that $(T,n)\leq (S,n)$ and for every comb $C\in \calC (T)$
there is $m$ such that $(T_{C(m)},m)\in D$.
This finishes the proof of Lemma 5 and Theorem 2.
\hspace{.25in} $\Box$

\section*{\bf Remarks}

In this paper we have only considered partitions of strongly increasing
sequences of reals and have shown that every such partition into 
an analytic and a co-analytic piece has a perfect homogeneous set.
A similar result can be obtained for any other type $\prec$ of increasing sequences.
All we have to do is modify the forcing notion ${\cal C}$ so that
the generic sequence produced has type $\prec$. 
Consequently, if ${\cal I}$ is a finite set of types of infinite increasing
sequences of reals for every analytic partitions of infinite increasing sequences
of reals we can find a perfect set $P$ such that for every  $\prec$ in ${\cal I}$
all the sequences from $P$ which have type $\prec$ have the same color.
On the other hand it is easy to see that if ${\cal I}$ is any infinite
set of types there is a partitions such that no perfect set is homogeneous
for all types in ${\cal I}$ simultaneously. Namely, choose for each
$s\in 2^{<\omega}$ a type $\prec _s$ in ${\cal I}$ such that the 
function which maps $s$ to $\prec _s$ is 1-1. Now given a sequence
$\{ x_n:n<\omega \}$ of type $\prec _s$ color it {\em red} if 
$\Delta (x_0,x_1)=s$ and {\em blue} otherwise. Let  now  $P$ be a perfect set 
and let  $T$ be  the tree of all finite initial segments of elements of $P$.
Then for any $s$ which is a splitting node of $T$ there are sequences
from $P$ of type $\prec _s$ which are colored by either color.

\vspace{.2in}

\indent  Universit\'e Paris VI, 

 Hebrew University, Jerusalem,  and 

 York University, Toronto

\end{document}